\documentclass[12pt]{article}

\usepackage{amsmath,theorem,amssymb,amsfonts,graphicx,color}

\newtheorem{theorem}{Theorem}[section]

\newtheorem{corollary}[theorem]{Corollary}

\newtheorem{lemma}[theorem]{Lemma}

\newtheorem{remark}[theorem]{Remark}

\newcommand{\dst}{\displaystyle}

\newcommand{\proof}{{\bf Proof} \hspace{2ex}}
\definecolor{red}{rgb}{1,0,0}
\definecolor{blue}{rgb}{0,0,1}

\begin{document}
\title{Properties of sub-matrices of Sylvester matrices and triangular toeplitz matrices}
\author{Yousong Luo, Robin Hill and Uwe Schwerdtfeger \\
{\small School of Mathematical and Geospatial Sciences, } \\ {\small RMIT University,
GPO Box 2476V} \\ {\small Melbourne, Vic. 3001, AUSTRALIA} \\
 email: {\small yluo@rmit.edu.au, r.hill@rmit.edu.au, u.schwerdtfeger@rmit.edu.au }}

\date{}
\maketitle

\begin{abstract}
In this note we discover and prove some interesting and important relations among sub-matrices of Sylvester matrices and triangular toeplitz matrices.  The main result is Hill's identity discovered by R. D. Hill which has an important application in optimal control problems.
\end{abstract}

\section{Introduction}

When studying the optimal state evolution of the dual state in a optimal control problem, R. Hill discovered an interesting relation (see Theorem \ref{robin}) among the sub-matrices of Sylvester matrices and triangular toeplitz matrices, see \cite{robin:1} and \cite{robin:2} for details.  If these relations holds then we can formulate the exact pattern how the modified states evolve.  In such a sense, the result here is not only an interesting result in linear algebra but also has a direct significant impact in control theory. 

We would like also to announce that we have an alternative proof for Theorem \ref{robin} using the tools given in \cite{Bezou:1} which is an entirely different approach.  

We formulate the problems first.  Define the following $m\hspace{-0.5mm} \times\hspace{-0.5mm} m$ lower and upper triangular matrices:
	\[ D_L :=  \left(
\begin{array}{ccccc}
  d_1&  0 & \cdots & \cdots & 0  \\
  d_2 & d_1  & 0 &\cdots & 0 \\
   \vdots & \ddots  & \ddots &\ddots & \vdots \\
    d_{m-1} & \ddots & \ddots &\ddots & 0  \\
      d_m & d_{m-1}  & \cdots &d_2 & d_1
\end{array}
\right)  \qquad
D_U :=  \left(
\begin{array}{ccccc}
  d_{m+1}&   d_m & \cdots & d_{3} &  d_{2}  \\
 0 &  d_{m+1}  & d_m  &\cdots &  d_{3}  \\
   \vdots & \ddots  & \ddots &\ddots & \vdots \\
   \vdots & \ddots & \ddots &\ddots &  d_m   \\
      0 & \cdots  & \cdots &0& d_{m+1}
\end{array}
\right)
\]
\[ N_L :=  \left(
\begin{array}{ccccc}
  n_1&  0 & \cdots & \cdots & 0  \\
  n_2 & n_1  & 0 &\cdots & 0 \\
   \vdots & \ddots  & \ddots &\ddots & \vdots \\
    n_{m-1} & \ddots & \ddots &\ddots & 0  \\
      n_m & n_{m-1}  & \cdots &n_2 & n_1
\end{array}
\right)  \qquad
N_U :=  \left(
\begin{array}{ccccc}
  n_{m+1}&   n_m & \cdots & n_{3} &  n_{2}  \\
 0 &  n_{m+1}  & n_m  &\cdots &  n_{3}  \\
   \vdots & \ddots  & \ddots &\ddots & \vdots \\
   \vdots & \ddots & \ddots &\ddots &  n_m   \\
      0 & \cdots  & \cdots &0& n_{m+1}
\end{array}
\right)
\]
Consider the Sylvester matrix
\[ S:=\left(\begin{array}{cc}D_L & N_L \\D_U &N_U \end{array}\right) \]
and the lower triangular matrix
\[ D:=\left(\begin{array}{cc}D_L & 0 \\D_U &D_L \end{array}\right). \]
The entries $d_1,d_2,\ldots, d_m, d_{m+1}$ and $n_1,n_2,\ldots, n_m, n_{m+1}$ are assumed to be nonzero real numbers such that both $S$ and $D$ are invertible.  Under such an assumption we define
\[ A:= D^{-1} \qquad
B:=S^{-1}. \]
If we use $A_T$ and $B_T$ to denote the matrices consisting of the first $m$ rows of $A$ and $B$, $A_B$ and $B_B$ the last $m$ rows of $A$ and $B$ respectively, then we can write
\[ A= \left(\begin{array}{c}A_T  \\A_B \end{array}\right) \qquad \mbox{and} \qquad
B= \left(\begin{array}{c}B_T  \\B_B \end{array}\right).\]
The $m\hspace{-0.5mm} \times\hspace{-0.5mm} m$ sub-matrices of $A_B$ consisting of the $m$ consecutive columns of it and starting from the $i$th column is denoted by $A_i$.  There are $m+1$ of them:
\begin{equation} \label{defa} A_1, A_2, \ldots,A_m , A_{m+1}.\end{equation}
Similarly, the sub-matrices of $B_B$ consisting of $m$ consecutive columns of it and starting from the $i$th column is denoted by $B_i$:
\begin{equation} \label{defb} B_1,B_2, \ldots,B_m,  B_{m+1}.\end{equation}

Our objective of this paper is to prove these relations, as well as discover and prove some other new relations among those sub-matrices.   The main result is the following Hill's identity.
\begin{theorem} \label{robin}  For $ 1 \leq i < j \leq m+1$ we have
\begin{equation}
A_iB_j=A_jB_i.
\end{equation}
\end{theorem}

The other results are
\begin{theorem} \label{relation1}
Assume that both $S$ and $D$ be invertible. Let $A_i$ and $B_j$ be the sub matrices defined in (\ref{defa}) and (\ref{defb}). Then, for all $i, j=1,\ldots m+1$, $A_i$ and $B_j$ are invertible and the following identities hold
\begin{equation} \label{non} {A_i}^{-1}A_j={B_i}^{-1}B_j \end{equation}
or equivalently
\begin{equation} \label{noi}A_j{B_j}^{-1}=A_i{B_i}^{-1}. \end{equation}
\end{theorem}
and
\begin{theorem}  \label{relation2} For $ 1 \leq i < j \leq m+1$ we have
\begin{equation}
B_i^{-1}B_j=B_jB_i^{-1}.
\end{equation}
\end{theorem}

As we can easily see that Theorem \ref{robin} is a consequence of the combination of Theorem \ref{relation1} and \ref{relation2}.

\section{Proofs of the results}
Now we introduce an $m\hspace{-0.5mm} \times\hspace{-0.5mm}3 m$ matrix
\begin{equation}
T:=
\begin{array}{c}
   \vspace{-2ex} \phantom{ (- D_U{D_L}^{-1} ,} \overbrace{\phantom{I_m, -D_L {D_U}^{-1}}}\\
  \vspace{-2ex} (\left. - D_U{D_L}^{-1} \; \right| \; I_m \; \left| \; -D_L {D_U}^{-1} \right.) \\
   \underbrace{\phantom{ (- D_U{D_L}^{-1} , I_m,}} \phantom{-D_L {D_U}^{-1}}
\end{array}
\end{equation}
where the symbol $|$ stands for an augmentation bar.  This matrix $T$ plays a very important role in the following argument through out this paper, so we call it ``kernel''. The $m\hspace{-0.5mm} \times\hspace{-0.5mm} 2m$ sub-matrices of $T$ consisting of the $2m$ consecutive columns of it and starting from the $i$th column is denoted by $T_i$ and we have $m+1$ such matrices:
\[ T_1, T_2, \ldots, T_m,  T_{m+1}.\]
Obviously $T_1=(- D_U{D_L}^{-1} , I_m)$ and $T_{m+1}=( I_m, -D_L {D_U}^{-1})$. Also, For each $i, j=1,2, \ldots, m+1$, the $m\hspace{-0.5mm} \times\hspace{-0.5mm} m$ sub-matrices of $T_i$ consisting of the $m$ consecutive columns of it and starting from the $j$th column is denoted by $T_{ij}$.

\begin{lemma} \label{lem1}If $\dst K=\left(\begin{array}{cc}D_L &0\\ D_U&D_L \\ 0 & D_U\end{array}\right)$, then
\begin{equation} \label{lm11}
TK =0.
\end{equation}
If $\dst D_l=\left(\begin{array}{c}D_L \\ D_U \end{array}\right)$, then for $i=1, 2, \ldots, m+1$ we have
\begin{equation} \label{lm12}
T_i D_l =0.
\end{equation}
\end{lemma}

\proof Obviously
\begin{eqnarray*}
TK &= & \left(\begin{array}{ccc} - D_U{D_L}^{-1} & I_m  & -D_L {D_U}^{-1}\end{array}\right) \left(\begin{array}{cc}D_L &0\\ D_U&D_L \\ 0 & D_U\end{array}\right) \\
&=&   \left(\begin{array}{cc} - D_U{D_L}^{-1}D_L+D_U & D_L- D_L{D_U}^{-1}D_U \end{array}\right) = \left(\begin{array}{cc} 0 & 0 \end{array}\right) .
\end{eqnarray*}
This immediately implies, by considering the first $m$ columns and the last $m$ columns of $TK$,  that
\begin{equation} \label{e1}
T_1 D_l = 0 \qquad \mbox{and} \qquad T_{m+1} D_l = 0 .
\end{equation}

For  $1<i< m+1$ let $K_i$ be the $m$ consecutive columns of $K$ starting from the $i$th column.  Then $K_i$ is in the form
\[ K_i = \left(\begin{array}{c}O_i \\ D_l \\ O_{m-i} \end{array}\right)\]
where $O_i$ is an $i \times m$ zero matrix and  $O_i$ is an $(m-1)i \times m$ zero matrix.
Therefore
\begin{equation} \label{e1}
T_i D_l =TK_i= 0.
\end{equation}
\null \hfill { QED}

{\bf Proof of Theorem \ref{relation1}} \hspace{2ex} We define
\begin{equation} \label{drdef}  D_r:=\left(
\begin{array}{c}
  0  \\
  D_{L} \end{array}
\right) \end{equation}
and hence
\[ D = \left( \begin{array}{cc}
  D_l & D_r  \end{array}
\right).\]
By Lemma (\ref{lem1}), $T_iD_l=0$.  Then, for $i, j = 1, \ldots, m+1$,  we have
\begin{equation} \label{keyeqn}
T_i= T_i D A = T_i \left(\begin{array}{cc} D_l & D_r \end{array}\right)A =
  \left(\begin{array}{cc} 0 & T_iD_r \end{array}\right) \left(\begin{array}{c} A_T \\ A_B \end{array}\right)=T_iD_rA_B
\end{equation}
which implies
\[ T_{ij}=T_iD_rA_j.\]
From the definition of $T$ we can see that $T_{m-i+2,i}=I$.  Then we have
\[I=T_{m-j+2}D_rA_j,\]
that is $A_j$ is invertible and
\begin{equation} \label{iaj}
{A_j}^{-1} =T_{m-j+2}D_r
\end{equation}
or
\begin{equation} \label{iaj2}
T_{i}D_r={(A_{m-i+2})}^{-1}.
\end{equation}
By substituting (\ref{iaj2}) into (\ref{keyeqn}) we obtain
\begin{equation} \label{keyrelationa}
T_{i}={(A_{m-i+2})}^{-1}A_B \qquad \mbox{or} \qquad {A_{i}}^{-1}A_B = T_{m-i+2}.
\end{equation}
This implies that
\begin{equation} \label{ainversea}
\qquad {A_{i}}^{-1}A_j = T_{m-i+2,j}.
\end{equation}

On the other hand we perform the same process to $B$ as follows.  We define
\begin{equation} \label{ndef}  N:=\left(
\begin{array}{c}
  N_L  \\
  N_{U} \end{array}
\right). \end{equation}

By Lemma (\ref{lem1}) we have, for $i, j = 1, \ldots, m+1$,
\begin{equation} \label{keyeqnb}
T_i = T_iSB = T_i \left(\begin{array}{cc} D_l & N \end{array}\right)B =
  \left(\begin{array}{cc} 0 & T_iN \end{array}\right) \left(\begin{array}{c} B_T \\ B_B \end{array}\right)=T_iNB_B
\end{equation}
which implies
\[ T_{ij}=T_iNB_j.\]
From the definition of $T$ we know that $T_{m-i+2,i}=I$.  Then we have
\[I=T_{m-j+2}NB_j,\]
that is
\begin{equation} \label{ibj}
T_{m-j+2}N={B_j}^{-1}
\end{equation}
or
\begin{equation} \label{ibj2}
T_{i}N={(B_{m-i+2})}^{-1}.
\end{equation}
By substituting (\ref{ibj2}) into (\ref{keyeqnb}) we obtain
\begin{equation} \label{keyrelationb}
T_{i}={(B_{m-i+2})}^{-1}B_B \qquad \mbox{or} \qquad {B_{i}}^{-1}B_B = T_{m-i+2}.
\end{equation}
This implies that
\begin{equation} \label{binverseb}
\qquad {B_{i}}^{-1}B_j = T_{m-i+2,j}.
\end{equation}
Equations (\ref{ainversea}) and (\ref{binverseb}) show that
\[ {A_{i}}^{-1}A_j={B_{i}}^{-1}B_j \]
for each $i, j = 1,2, \ldots, m+1$.  This completes the proof. \null \hfill { QED}

\begin{corollary}  \label{coro1} We define
\begin{equation} \label{mdef}  M:=\left(
\begin{array}{cc}
 M_{1} & M_2 \end{array}
\right) = \left(
\begin{array}{cc}
  N_L & 0 \\
  N_{U}& N_L\\
  0 & N_{U} \end{array}
\right).
\end{equation}
Let $ H = T M $ and $H_i$ be the sub-matrix of $H$ consisting the $m$ consecutive columns of $H$ starting from the $i$th column.  Then
\[ H_i = {(B_{m-i+2})}^{-1} \qquad \mbox{or} \qquad  H_{m-i+2}= {B_{i}}^{-1}.\]
\end{corollary}

\proof Consider
\begin{equation} \label{h1}
H=TM =T\left(
\begin{array}{cc}
  N_L & 0 \\
  N_{U}& N_L\\
  0 & N_{U} \end{array}
\right) =\left(  \begin{array}{cc}  T_1N & T_{m+1} N \end{array}
\right).
\end{equation}
This gives immediately
\begin{equation} \label{h2}
H_1=T_1 N \qquad \mbox{and} \qquad  H_{m+1}=T_{m+1} N.
\end{equation}
Equations (\ref{ibj2}) then implies $H_1={(B_{m+1})}^{-1}$ and $H_{m+1}={B_1}^{-1}$.
For  $1<i< m+1$ let $M_i$ be the sub-matrix of $M$ consisting the $m$ consecutive columns of $M$ starting from the $i$th column.  Then $M_i$ is in the form
\[ M_i = \left(\begin{array}{c}O_i \\ N \\ O_{m-i} \end{array}\right)\]
where $O_i$ is an $i \times m$ zero matrix and  $O_i$ is an $(m-1)i \times m$ zero matrix.
Therefore
\begin{equation} \label{e1}
H_i=TM_i= T_iN.
\end{equation}
Again, equations (\ref{ibj2}) shows $H_i={(B_{m-i+2})}^{-1}$.
\null \hfill { QED}

\begin{remark}
This theorem reveals two remarkable features of $A_i$'s and $B_i$'s. First, equation (\ref{noi}) demonstrates the invariance of $A_i{B_i}^{-1}$ with respect to $i$.  More precisely we have
    \[ A_i{B_i}^{-1} = A_B N.\]
Secondly, equation (\ref{non}) shows that ${B_i}^{-1}B_j$ is independent of $n_h$'s which are the elements defining $S$.  This is quite significant as $B_i$'s are sub-matrices of $B$, which is the inverse of $S$ and therefore depends on $n_h$'s.
\end{remark}

\begin{remark}
The proof of this theorem also demonstrates an interesting feature of those $A_i$'s and $B_i$'s.  By the definition of $T$ we can see that, for $i, j =1, 2, \ldots, m+1$ and $1 \leq k \leq \max \{ m-i+1, j\}$ we have
    \[ T_{i+k,j-k} = T_{i,j}.\]
This, together with (\ref{ainversea}) and (\ref{binverseb}), shows that
\begin{equation} \label{rmk2} {A_{i}}^{-1}A_j = {(A_{i+k})}^{-1}A_{j+k} \qquad \mbox{and} \qquad  {B_{i}}^{-1}B_j = {(B_{i+k})}^{-1}B_{j+k}\end{equation}
for such $k$'s that the right hand sides of the above equations are defined.
For example,
\[{B_{1}}^{-1}B_2 ={B_{2}}^{-1}B_3= \cdots = {B_{m}}^{-1}B_{m+1}.\]
\end{remark}

{\bf Proof of Theorem \ref{relation2}} \hspace{2ex} It is well known that $B$ can be represented by
\begin{equation} \label{bzform}
B=  \left(\begin{array}{cc} N_UB_z & -N_LB_z \\ -D_UB_z & D_LB_z\end{array}\right)
\end{equation}
where $B_z={B_T(D,N)}^{-1}$ where $B_T(D,N)$ is the Bezoutian matrix generated by $D$ and $N$ in the following manner:
\begin{equation} \label{bz}
 B_T(D,N) = D_LN_U-N_LD_U =N_UD_L-D_UN_L.
\end{equation}
For detailed properties of Bezoutian matrices we refer to the comprehensive article \cite{Bezou:1}. Using this representation we have $B_1= -D_UB_z$ and $B_{m+1}=D_LB_z$.

Now, by Corollary \ref{coro1}, we have
\begin{eqnarray*}
B_1H
&=&  B_1\left(  \begin{array}{cc}  {(B_{m+1})}^{-1}  &
 {B_{1}}^{-1}  \end{array} \right) = \left(  \begin{array}{cc}  B_1{(B_{m+1})}^{-1}  &
I \end{array} \right)\\
  &=&  \left(  \begin{array}{cc}   -D_UB_z ({B_z}^{-1} {D_L}^{-1}) \;  &
  \;  I \end{array} \right) =\left(  \begin{array}{cc}  -D_UD_L \;  &
  \;  I \end{array} \right) \\
   &=& T_1
\end{eqnarray*}
and hence
\[ B_1{B_i}^{-1}=B_1H_{m-i+2}=T_{1, m-i+2}.\]
This, together with equation (\ref{keyrelationb}), implies
\[  B_1{B_i}^{-1}={(B_{m+1})}^{-1}B_{m-i+2}.\]
Putting $k=m-i+1$ in (\ref{rmk2}) gives
\[ {B_{i}}^{-1} B_1={(B_{i+k})}^{-1}B_{1+k}={(B_{m+1})}^{-1}B_{m-i+2}.\]
Therefore $B_1B_i^{-1} =B_i^{-1}B_1$ for each $i=1,2, \ldots, m+1$.

Similarly
\begin{eqnarray*}
B_{m+1}H
&= & B_{m+1} \left(  \begin{array}{cc}  {(B_{m+1})}^{-1}  &
 {B_{1}}^{-1}  \end{array} \right) = \left(  \begin{array}{cc}  I  &
 B_{m+1}{B_1}^{-1}  \end{array} \right) \\
  & =  &  \left(  \begin{array}{cc}  I \;  &   \;   D_LB_z (-{B_z}^{-1} {D_U}^{-1})  \end{array} \right) =\left(  \begin{array}{cc} I \;  &
  \; -D_L{D_U}^{-1} \end{array} \right) \\
   &=& T_{m+1}.
\end{eqnarray*}

This, together with equation (\ref{keyrelationb}), proves
\[ B_{m+1} {B_i}^{-1} = T_{m+1, m+2-i} ={B_1}^{-1}B_{m+2-i}.\]
Equation (\ref{rmk2}) with $k=i-1$ gives
\[{B_1}^{-1}B_{m+2-i}={(B_{1+i-1})}^{-1}B_{m+2-i+i-1}={B_{i}}^{-1}B_{m+1},\]
and hence $B_{m+1}B_i^{-1} =B_i^{-1}B_{m+1}$ for each $i=1,2, \ldots, m+1$.  This is equivalent to
\begin{equation} \label{bibyblast}
 B_i{(B_{m+1})}^{-1} ={(B_{m+1})}^{-1}B_i.
\end{equation}

Now for $1<i<m+1$, by equation (\ref{binverseb})
\begin{eqnarray*}
B_iH
&= & B_i \left(  \begin{array}{cc}  {(B_{m+1})}^{-1}  &
 {B_{1}}^{-1}  \end{array} \right) = \left(  \begin{array}{cc}  B_i{(B_{m+1})}^{-1} &
 B_i{B_1}^{-1}  \end{array} \right) \\
  & =  &  \left(  \begin{array}{cc}  {(B_{m+1})}^{-1}B_i \;  &   \;   {B_1}^{-1}B_i   \end{array} \right) \\
   &=& \left(  \begin{array}{cc}  T_{1,i} \;  &    T_{m+1,i}  \end{array} \right).
\end{eqnarray*}
Let $t_j$ denote the $j$th column of $T$.  The observation
\begin{equation}
T =
\begin{array}{ccc}
    (t_1, \ldots, t_{i-1},& \underbrace{\overbrace{t_i, \ldots, t_{m+i-1}}^{ T_{1,i}}, \overbrace{t_{m+i},  \ldots,t_{2m+i-1}}^{ T_{m+1,i}}}_{T_i}, & t_{2m+i}, \ldots, t_{3m} )
\end{array}
\end{equation}
shows that
\[\left(  \begin{array}{cc}  T_{1,i} \;  &    T_{m+1,i}  \end{array} \right) = T_i,\]
and hence
\begin{equation} \label{important2}
B_iH = T_i.
\end{equation}
From this we obtain  $B_j{B_i}^{-1} ={B_i}^{-1}B_j$. \hfill QED
\begin{corollary} For $i, j = 1, 2, \ldots, m+1$ we have
\begin{equation} \label{bibyblast}
 B_iB_{j} =B_{j}B_i,
\end{equation}
and, for all $l$ such that both $B_{i+l}$ and $B_{j-l}$ are meaningful,
\begin{equation} \label{bibyblast}
 B_iB_{j} =B_{i+l}B_{j-l}.
\end{equation}
\end{corollary}
\proof The second equation follows from (\ref{rmk2}) by putting $k=i-j+l$:
\[ B_{j-l}{B_{j}}^{-1}={B_{j}}^{-1}B_{j-l}={(B_{j+k})}^{-1}B_{j-l+k}={(B_{i+l})}^{-1}B_{i}.\]


\begin{thebibliography}{}
%
%
\bibitem{Bezou:1} Georg Heinig and Karla Rost,
{\em Introduction to Bezoutians},
Advances and Applications, Vol. 199, 25 - 118, (2010)

\bibitem{robin:1} Robin Hill, Uwe Schwerdtfeger and Michael Baake, {\em Dynamic programming and duality applied to an optimal control problem},
Proceeding of Australian Control Conference, to appear, (2011)

\bibitem{robin:2} Robin D. Hill, {\em Dual periodicity in $l_1$-norm minimisation problems},
Systems \& Control Letters, 57, 489 - 496, (2008)

\end{thebibliography}
\end{document}